\documentclass[12pt]{article}
\usepackage{amssymb,amsmath, amsthm, amscd,ifthen}
\usepackage[T2A]{fontenc}
\usepackage[cp1251]{inputenc}
\usepackage[english,russian]{babel}
\usepackage[dvips]{graphicx}
\usepackage[hyper]{amsbib}
\usepackage{indentfirst}
\usepackage{enumerate}

\usepackage{amsmath}
\usepackage{amssymb}

\newtheorem{theorem}{Theorem}
\newtheorem{claim}{Claim}

\newtheorem{lemma}{Lemma}

\newenvironment{proof1}{\noindent{\it Proof.\,}}{\hfill$\Box$}

\begin{document}

\title{Equitable colorings of hypergraphs with few edges}

\author{Margarita B.\,~Akhmejanova\footnote{Moscow Institute of Physics and Technology, Laboratory of Advanced Combinatorics and Network Applications, 141700, Institutskiy per. 9, Dolgoprudny, Moscow Region, Russia;}, Dmitry A.\,~Shabanov\footnote{Moscow Institute of Physics and Technology, Laboratory of Advanced Combinatorics and Network Applications, 141700, Institutskiy per. 9, Dolgoprudny, Moscow Region, Russia; 
Lomonosov Moscow State University, Faculty of Mechanics and Mathematics, Department of Probability Theory, 119991, Leninskie gory, 1, Moscow, Russia; 
National Research University Higher School of Economics (HSE), Faculty of Computer Science, 101000, Myasnitskaya Str. 20, Moscow, Russia. E-mail: dmitry.shabanov@phystech.edu}}

\date{}
\maketitle

\begin{center}
\textbf{Abstract}
\end{center}

The paper deals with an extremal problem concerning equitable colorings of uniform hyper\-graph. Recall that a vertex coloring of a hypergraph $H$ is called proper if there are no monochro-matic edges under this coloring. A hypergraph is said to be equitably $r$-colorable if there is a proper coloring with $r$ colors such that the sizes of any two color classes differ by at most one. In the present paper we prove that if the number of edges $|E(H)|\leq 0.01\left(\frac{n}{\ln n}\right)^{\frac {r-1}{r}}r^{n-1}$ then the hypergraph $H$ is equitably $r$-colorable provided $r<\sqrt[5]{\ln n}$.

\section{Introduction}

The present paper deals with the well-known problem concerning colorings of hypergraphs. Let us start with recalling some definitions.

\subsection{Main definitions and new result}
Let $H=(V,E)$ be a hypergraph. A coloring of the vertex set $V$ of a hypergraph $H=(V,E)$ is called \emph{proper} if none of the edges in $E$ is monochromatic under this coloring. A hypergraph is said to be $r$-\emph{colorable} if there exists a proper coloring with $r$ colors for it. A coloring of the hypergraph vertices is said to be \emph{equitable} if it is proper and the sizes of any two of color classes differ by at most one. The last means that the set of vertices $V$ can be divided not just into $r$ independent sets, but into $r$ independent sets of almost the same size.

\bigskip
The main result of the paper provides a sufficient condition for admitting an equitable $r$-coloring by an $n$-uniform hypergraph as the restriction on the number of edges.

\begin{theorem}\label{thm:main}
	For large enough $n$, if $r<\sqrt[5]{\ln n}$ then any $n$-uniform hypergraph $H=(V,E)$ with
	\begin{equation*}
	|E|\leq 0.01\left(\frac{n}{\ln n}\right)^{\frac {r-1}{r}}r^{n-1}
	\end{equation*}
	is equitably $r$-colorable provided $|V|$ is divisible by $r$.
\end{theorem}

\subsection{Related work and method}

Let us begin with well-known extremal problem concerning hypergraphs colorings. This whole thing started when Erd\H{o}s and Hajnal proposed to find the value $m(n,r)$ which is equal to the smallest number of edges in an $n$-uniform none-$r$-colorable hypergraph. We will only note that for small $r$ in comparison with $n$, the best current bounds are the following:
\begin{equation}\label{bound_CherkKozik}
  c_1\left(\frac{n}{\ln n}\right)^{\frac{r-1}r}r^{n-1}\le m(n,r)\le c_2 n^2r^n\ln r,
\end{equation}
were $c_1,c_2>0$ are some absolute constants. The lower bound was proved by Kozik and Cherkashin\cite{CherkKozik}, the proof of upper bound is that of Erd\H{o}s (see the proof in Kostochka's paper \cite{Kost}). In fact, the result of Kozik and Cherkashin remains true for large $r$ compared to $n$, but in this case, Akolzin and Shabanov proved more stronger estimates \cite{AkolzinShabanov}. Reader can find more results on the related problems in \cite{Cherkashin}--\cite{AkhShab}. Recent advances concerning colorings of random hypergraphs were obtained in \cite{DFG}--\cite{KKSh}.

\bigskip
Now we will move on to history of equitable colorings. A significant interest in equitable colorings  historically arose in connection with the celebrated theorem of Hajnal and Szemer\'edi \cite{HajSem} of 1970, which verified the following conjecture of Erd\H{o}s \cite{Erd2}: \emph{every graph $G$ with maximal vertex degree $\Delta(G)$ admits not only a proper but also an equitable coloring with $\Delta+1$ colors}.

It is worth noting that the original proof of Hajnal and Szemer\'edi \cite{HajSem} was long enough and complicated. In 2008 Kierstead and Kostochka \cite{KierKost} presented another proof, that was short and transparent. Moreover, they refined this result to the case of maximum edge degree \cite{KierKost2} and together with  Mydlarz and Szemer\'edi found a fast algorithm for obtaining an equitable coloring \cite{KierKostMyd}.

The generalization of Erd\H{o}s's conjecture to uniform hypergraphs was obtained by Lu and Sz\'ekely. Namely, they proved that:\emph{ if $|V|$ is divisible by $r$ and the maximum vertex degree $\Delta(H)\leq\ r^{n-1}/(2en),$
	then  $H$ admits an equitable $r$-coloring} \cite{LuSz}.

\bigskip
The main aim of the present paper is to give the sufficient condition for an equitable $r$-colorability of a hypergraph as a restriction on the number of edges in the general class of $n$-uniform hypergraphs. In other words, we study the analog of the Erd\H{o}s--Hajnal problem for equitable colorings \cite{ErdHaj}. We also stress, that the previous best known bound that guarantees the existence of an equitable $r$-coloring of an $n$-uniform hypergraph was just $r^{n-1}$.

The proof of Theorem \ref{thm:main} is based on the idea that there is a coloring in which there are few specific configurations, so-called ordered $k$-chains, where $k\in\{1,2,...,r\}$. The idea of considering ordered $r$-chains was first conceived and proved by Pluh\H{a}r \cite{Pluhar}: \emph{an arbitrary hypergraph is $r$-colorable if and only if there exists order on $V$ without ordered $r$-chains}. The similar idea was used by Cherkashin and Kozik \cite{CherkKozik}, who get the best known result for the case when $n$ is big enough in comparison with $r$.

Let's get back to Theorem \ref{thm:main}. The crucial moment is that:\emph{ Pluh\H{a}r criterion
does not provide the information on the cardinalities of color classes, while consideration of ordered chains of all sizes does provide it.}

\bigskip
The rest of the paper is organized as follows. In Section 2 we show that the theorem holds for hypergraphs with few vertices. Section 3 is devoted to the random algorithm, which with positive probability creates proper, but not necessarily equitable $r$-coloring of $H$. In Section 4, by the recoloring of few vertices, we will ensure the equalities of sizes of the color classes keeping the lack of the monochromatic edges.

\section{Hypergraphs with few vertices}

Assume that $H=(V,E)$ is a hypergraph from Theorem \ref{thm:main}. Let $m=|V|$ denote its number of vertices. The aim of this paragraph is to show that the following conditions:
	\begin{equation}\label{vert_bound}
	|E|<0.01\left(\frac{n}{\ln n }\right)^{\frac{r-1}{r}} r^{n-1},~~~|V|<\frac{n^2(r-1)}{2\ln n},
	\end{equation}
imply that the hypergraph $H$ admits an equitable coloring with $r$ colors. For this purpose, it is sufficient to consider a random balanced coloring with $r$ colors. A coloring is said to be balanced if the sizes of its color classes are equal, i.e. it is a partition of $V$ into $r$ parts with the same size: $V=K_{1}\sqcup K_{1}...\sqcup K_{r-1}, |K_{1}|=|K_{2}|=...=|K_{r}|=m/r$.

\bigskip
Let $C$ be a random balanced coloring of $V$ with $r$ colors. Then for any edge $A\in E$,
$$
  \Prob(\text{$A$ is monochromatic under $C$})=\frac{r\binom{m-n}{m/r-n}}{\binom{m}{m/r}}.
$$
Thus,
\begin{align*}
	\Prob(&\text{there exists a monochromatic edge under $C$})\leq\\
&=|E|\frac{r\binom{m-n}{m/r-n}}{\binom{m}{m/r}}=|E|\frac{r^{-n+1}(1-r/m)...(1-r(n-1)/m)}{(1-1/m)...(1-(n-1)/m)}\leq\\
	&\leq 0.01\left(\frac{n}{\ln n}\right)^{\frac{r-1}{r}}\frac{\exp(\ln(1-r/m)+...+\ln(1-r(n-1)/m))}{\exp(\ln(1-1/m)+...+\ln(1-(n-1)/m))}=\\
	&= 0.01\left(\frac{n}{\ln n}\right)^{\frac{r-1}{r}}\prod_{x=1}^{n-1}\exp\left(\ln(1-rx/m)-\ln(1-x/m)\right)
	\end{align*}
	The Taylor expansions of the logarithmic functions imply that
$$
    \ln(1-rx/m)-\ln(1-x/m)<-x(r-1)/m
$$
for $\frac{rx}m\in(0,1)$. Summarizing the degrees in the product of exponents, we finally obtain:
	\begin{align*}
	\Prob(&\text{there exists a monochromatic edge under $C$})\leq \\
    &\leq 0.01\left(\frac{n}{\ln n}\right)^{\frac{r-1}{r}}e^{-n(n-1)(r-1)/2m}\leq\\
    &\leq 0.01 n^{\frac{r-1}{r}-\frac{n(n-1)}{n^2}}\cdot (\ln n)^{-\frac{r-1}{r}}=\\
    &= 0.01 n^{\frac {-1}{r}+\frac 1{n}}\cdot (\ln n)^{-\frac{r-1}{r}}<1.
	\end{align*}
Hence, with positive probability the random balanced coloring $C$ is equitable. It remains to consider hypergraphs with the number of vertices greater than $n^2(r-1)/(2\ln n)$.

\section{Algorithm 1: construction of a proper coloring}

Assume now that the hypergraph $H$ satisfies the conditions of Theorem \ref{thm:main} and has large number of vertices:
\begin{equation}\label{vert_bound_new}
m=|V|\ge\frac{n^2(r-1)}{2\ln n}.
\end{equation}
We denote the set of colors by $\{1,\ldots,r\}.$

For every vertex $v\in V$, let $\sigma(v)$ be an independent random variable with uniform distribution on $[0,1)$. The value $\sigma(v)$ is called \emph{the weight} of the vertex $v$. With probability 1 the mapping $\sigma:V\to [0,1)$ is an injection. So, $\sigma$ induces a random ordering on $V$, i.e. $(v_{1},v_{2},...,v_{m})$ are the vertices of $H$ written in the order $\sigma$ if
$\sigma(v_{1})<\sigma(v_{2})<\ldots<\sigma(v_{m})$. The Algorithm 1 will be parametrized by the value
\begin{equation}\label{choice_p}
p=\left(\frac{r-1}{r}\right)\frac{\ln(\frac{n}{\ln n})}{n}.
\end{equation}

\begin{figure}[h]
	\begin{minipage}[h]{1\linewidth}
		\center{\includegraphics[width=0.8\linewidth]{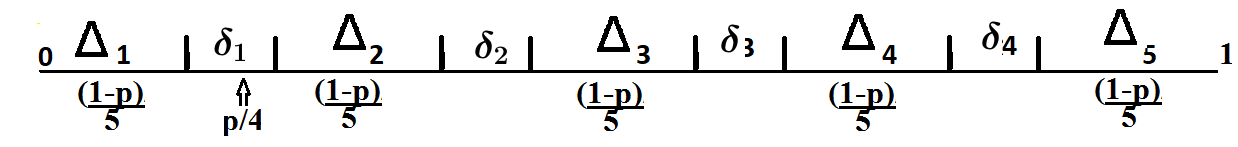}}
	\end{minipage}
	\selectlanguage{english}
\caption{Partition of $[0,1)$ into $\Delta_1,\delta_1,\Delta_2,\delta_2,\ldots,\Delta_5$ for $r=5$.}
\label{picture}
\end{figure}

\bigskip
We divide the unit interval $[0,1)$ into subintervals $\Delta_1,\delta_1,\Delta_2,\delta_2,\ldots,\Delta_r$ as on the Figure \ref{picture}, i.e.
$$
  \Delta_i=\left[(i-1)\left(\frac {1-p}r+\frac p{r-1}\right),i\cdot\frac {1-p}r+(i-1)\cdot\frac p{r-1}\right),\;i=1,\ldots,r;
$$
$$
  \delta_i=\left[i\cdot\frac {1-p}r+(i-1)\cdot\frac p{r-1},i\left(\frac {1-p}r+\frac p{r-1}\right)\right),\;i=1,\ldots,r-1.
$$
The length of each large subinterval $\Delta_i$ is equal to $\left(\frac{1-p}{r}\right)$, and every small subinterval $\delta_i$ has length equal to $\left(\frac{p}{r-1}\right)$.
A vertex $v$ is said to belong to a subinterval $[c,d)$, if its weight $\sigma(v)\in[c,d)$.

\bigskip

\bigskip
We color the vertices of hypergraph $H$ according to the following Algorithm 1, which consists of two stages.
\begin{itemize}
\item  First, each $v\in\Delta_{i}$ is colored with color $i$, $i=1,\ldots,r$
\item  Then, moving with growth of weight, we color a vertex $v\in\delta_{i}$ with color $i$ if such assignment does not create a monochromatic edge in the current coloring. Otherwise we color $v$ with color $i+1$.
\end{itemize}

Let $C^0$ denote the random coloring obtained after the consideration of all the vertices.

\subsection{Analysis of Algorithm 1}

Suppose that Algorithm 1 fails to produce a proper coloring and there is a monochromatic edge $A$ in the initial coloring $C^0$. Let $i\in\{1,2,...,r\}$ denote the color of $A$ and let $v_A$ be the first vertex of $A$. We note that $v_A$ could receive color $i$ only in two cases: either $v_A\in\Delta_i$ or $v_A\in\delta_{i-1}$. In the second case there exists an edge $B$, such that $v_A$ is the last vertex of $B$ and the remaining vertices of $B$ were colored with color $(i-1)$. In this situation we say that the pair $(A,B)$ is \emph{conflicting}.

For the first vertex $v_B$ of the edge $B$ we also have an alternative: either $v_B\in\Delta_{i-1}$ or $v_B\in\delta_{i-2}$ and there exists an edge $C$ such that $v_B$ is the last vertex of $C$ and the remaining vertices of $C$ are colored with color $(i-2)$. Repeating the above arguments, we obtain a construction, called \emph{an ordered $k$-chain for color $i$}. It is an edge sequence $H'=(C_1,...,C_{k-1}=B,C_k=A)$ such that the first vertex of $C_1$ belongs to the subinterval $\Delta_{i-k+1}$, and for every  $j=2,...,k$, pair $(C_{j-1}, C_j)$ is conflicting.

\bigskip
Let us make some important notes concerning ordered chains.
\begin{enumerate}
  \item The case of ordered $1$-chain corresponds to the case when $v_A\in\Delta_i$.
  \item The last vertex of the edge $C_k=A$ belongs to the subinterval $\Delta_i$ (otherwise we should prefer the color $i+1$).
  \item Every pair $(C_{j-1}, C_j)$ has exactly one vertex in common and this vertex belongs to $\delta_{j-1}$.
\end{enumerate}

Summarizing the above, we can say that
\begin{claim}\label{claim1}
 If for injective $\sigma:V\to[0;1)$ and for each color $i\in\{1,2,...,r\}$, there are no ordered chains then Algorithm 1 produces a proper coloring.
\end{claim}

\subsection{Auxiliary claims concerning ordered $k$-chains}

\begin{lemma}\label{lemma1}
The number of configurations in the hypergraph $H$ that can form an ordered $k$-chain is at most $2\binom{|E|}{k}$.
\end{lemma}

\begin{proof1}
Let us take an arbitrary unordered family of $k$ edges of $H$, $A_1,\ldots,A_k$. This can be made by at most $\binom{|E|}{k}$ ways. If this set can form a chain, say, $(A_1,\ldots,A_k)$, then any successive two edges have exactly one common vertex and the remaining pairs do not intersect. So, there is only one more possible chain $(A_k,\ldots,A_1)$.
\end{proof1}

\begin{lemma}\label{lemma_prob1}
Suppose $k\geq 1$ and let $H'=(C_1,\ldots,C_k)$ be an ordered $k$-tuple of edges in the hypergraph $H$. Then the probability that $H'$ forms an ordered $k$-chain for color $i$ does not exceed $$2\left(\frac{\ln n}{n}\right)^{\frac{k(r-1)}{r}} r^{-(n-1)k-1}.$$
\end{lemma}

\begin{figure}[h]
	\begin{minipage}[h]{1\linewidth}
\center{\includegraphics[width=0.8\linewidth]{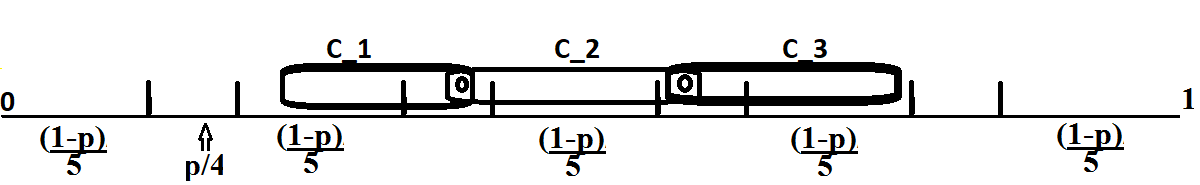}}
	\end{minipage}
	\selectlanguage{english}
	\caption{Ordered $3$-chain for color $4$}
	\label{picture2}
\end{figure}
\begin{proof1}
Suppose $k\geq 2$. Let $v_j=C_j\cap C_{j+1}$ denote a common vertex of $C_j$ and $C_{j+1}$, $j=1,\ldots,k-1$, i.e. the vertex $v_j$ is the last vertex of the edge $ C_j$ and is the first vertex of the edge $C_{j+1}$. Obviously, $v_j\in\delta_{i-k+j}$, $j=1,\ldots,k-1$. Let us also denote $\delta_j=[\alpha_j,\beta_j)$, $\beta_j-\alpha_j=p/(r-1)$.

So, for given weights $\sigma(v_j)=x_j$, $j=1,\ldots,k-1$, the event that $H'$ forms an ordered $k$-chain for color $i$, can be described as follows:
\begin{itemize}
  \item every vertex $w\in C_1\setminus\{v_1\}$ should belong to the subinterval $\Delta_{i-k+1}\sqcup[\alpha_{i-k+1},x_{1})$;
  \item every vertex $w\in C_j\setminus(\{v_j,v_{j-1}\})$, $j=2,\ldots,k-1$, should belong to the subinterval $[x_{j-1},\beta_{i-k+j-1})\sqcup\Delta_{i-k+j}\sqcup[\alpha_{i-k+j},x_{j})$;
  \item every vertex $w\in C_k\setminus\{v_{k-1}\}$ should belong to the subinterval $[x_{k-1},\beta_{i-1})\sqcup\Delta_{i}$.
\end{itemize}
The weights are independent, so denoting $y_j=x_j-\alpha_{i-k+j}$, we obtain the following estimate for the conditional probability:
\begin{align*}
&\left(\frac{1-p}{r}+y_1\right)^{n-1}\times\left(\frac{1-p}{r}+y_2+\frac{p}{r-1}-y_1\right)^{n-2}\times\ldots\times\\
&\times\left(\frac{1-p}{r}+y_{k-1}+\frac{p}{r-1}-y_{k-2}\right)^{n-2}\times\left(\frac{1-p}{r}+\frac{p}{r-1}-y_{k-1}\right)^{n-1}\leq\\
&\leq\text{|take out the factor $r^{-k(n-2)-2}$ and use the estimate $(1+y)^s\leq\exp\{ys\}$}|\leq\\
&\leq r^{-k(n-2)-2}\exp\{(n-1)\left(-p+ry_1\right)+(n-2)\left(-p+ry_2-ry_1+pr/(r-1)\right)+\ldots+\\
&+(n-2)\left(-p+ry_{k-1}-ry_{k-2}+pr/(r-1)\right)+(n-1)\left(-p+pr/(r-1)-ry_{k-1}\right)\}\leq\\
&\leq\text{|since $0\leq y_j\leq p/(r-1)$|}\leq\\
&\leq r^{-k(n-2)-2}e^{(n-2)\left(-pk+\frac{(k-1)pr}{r-1}\right)}\cdot e^{\frac{2p}{r-1}}.
\end{align*}
To obtain the final estimate, we have to integrate over $(y_1,\ldots\,y_{k-1})\in [0,p/(r-1)]^{k-1}$ and substitute $p$ from \eqref{choice_p}. Thus, the probability under the consideration does not exceed
\begin{align*}
\left(\frac{p}{r-1}\right)^{k-1}&r^{-k(n-2)-2}e^{(n-2)\left(-pk+\frac{(k-1)pr}{r-1}\right)}\cdot e^{\frac{2p}{r-1}}\leq\\
&\leq\left|\text{since $e^{2p/(r-1)}\leq 2$ and $(-pk+(k-1)pr/(r-1))=\left(-1+\frac{k}{r}\right)\frac{\ln\left(\frac{n}{\ln n}\right)}{n}$}\right|\leq\\
&\leq r^{-k(n-1)-1}\left(\frac{\ln n}{n}\right)^{k-1}\cdot e^{\left(-1+\frac{k}{r}\right)\ln\left(\frac{n}{\ln n}\right)}\cdot 2\leq\\ &\leq 2\left(\frac{\ln n}{n}\right)^{\frac{k(r-1)}{r}} r^{-(n-1)k-1}.
\end{align*}
If $k=1$ then every vertex of $C_1$ should belong to $\Delta_i$, so the probability does not exceed
$$
  \left(\frac{1-p}{r}\right)^n\leq\left(\frac{\ln n}{n}\right)^{\frac{r-1}{r}} r^{-n}.
$$
\end{proof1}

\subsection{Outcome of Algorithm 1}

Lemmas \ref{lemma1} and \ref{lemma_prob1} immediately imply the following statement.

\begin{lemma}\label{lemma_prob3} The probability that a monochromatic edge occurs in the random coloring
$C^0$ does not exceed $0.04e$.
\end{lemma}
\begin{proof1} Lemma \ref{lemma1} estimates the number of possible configurations that can form an ordered
$k$-chain and Lemma \ref{lemma_prob1} does the same for the probability that a $k$-chain happens for a given color $i$. Bringing them together, we obtain that the required probability does not exceed
$$
  \sum_{k=1}^r 2\binom{|E|}{k}\cdot 2r\left(\frac{\ln n}{n}\right)^{\frac{k(r-1)}{r}}\left(\frac 1r\right)^{(n-1)k+1}\leq 4\sum_{k=1}^\infty\frac{(0.01)^k}{k!}\leq 0.04 \cdot\sum_{k=1}^\infty\frac{1}{k!}<0.04e.
$$
\end{proof1}

\subsection{Some auxiliary claims concerning color classes}

Let us return to the analysis of Algorithm 1:

\begin{itemize}
\item  First, each $v\in\Delta_{i}$ is colored with color $i$, $i=1,\ldots,r$
\item  Then, moving with growth of weight, we color a vertex $v\in\delta_{i}$ with color $i$ if such assignment does not create a monochromatic edge in the current coloring. Otherwise we color $v$ with color $i+1$.
\end{itemize}

Let $X(i)$ denote the number of vertices $v\in\delta_{i}$ which have been colored with color $i+1$ during the evaluation of Algorithm 1. In other words, $X(i)$ is equal to the number of vertices $v\in\delta_{i}$ which Algorithm 1 could not color with $i$. The following lemma provides an estimate of the expected values of $X(i)$ in terms of number of edges $|E|$.
\begin{lemma}\label{lemma4}
	Suppose that number of edges $|E|$ of $H$ is at most  $0.01\left(\frac{n}{\ln n}\right)^{\frac{r-1}{r}} r^{n-1}$. Then
\begin{enumerate}[(i)]
\item the expected value of $X(i)$ does not exceed $0.04e\cdot \frac{n}{r\ln n}$, $i=1,\ldots,r-1$;
\item with probability at least $(1-0.04e)$ every $X(i)$ does not exceed $\frac{n}{\ln n}$.
\end{enumerate}
\end{lemma}
\begin{proof1}
Suppose that $v\in\delta_i$ was colored with color $i+1$ during the evaluation of Algorithm 1. Since the color $i$ was not allowed, there exists an edge $A$ such that if $v$ is colored with color $i$ then $A$ would become monochromatic of color $i$. This means that all the vertices in $A\setminus\{v\}$ have smaller weights than $v$ and have been colored with $i$ during Algorithm 1. Now we can consider the first vertex of $A$. The same argument as in the proof of Claim \ref{claim1} provides a chain $(C_1,\ldots,C_k)$ such that
\begin{enumerate}
  \item $C_k=A$.
  \item Every pair $(C_j,C_{j+1})$, $j=1,\ldots,k-2$, is a conflicting pair. A vertex $v_j=C_j\cap C_{j+1}$ belongs to $\delta_{i-k+j}$.
  \item The first vertex of $C_1$ belongs to $\Delta_{i-k+1}$.
  \item A pair $(C_{k-1},C_k)$ has one common vertex $v_{k-1}\in\delta_{i-1}$.
  \item The last vertex $v$ of the edge $C_k$ belongs to $\delta_i$.
\end{enumerate}
Roughly speaking, the last property is the only difference with ordered $k$-chain for color $i$. We will say that in the situation described above the set of edges $(C_1,\ldots,C_k)$ forms \emph{an improper k-chain for color $i$}.

\bigskip
Estimation of the probability almost repeats the argument in Lemma \ref{lemma_prob1}. Again we use the notation $\delta_j=[\alpha_j,\beta_j)$, $\beta_j-\alpha_j=p/(r-1)$. For given weights $\sigma(v_j)=x_j$, $j=1,\ldots,k-1$, the event that $(C_1,\ldots,C_k)$ forms an improper ordered $k$-chain for color $i$, can be described as follows:
\begin{itemize}
  \item every vertex $w\in C_1\setminus\{v_1\}$ should belong to the subinterval $\Delta_{i-k+1}\sqcup[\alpha_{i-k+1},x_{1})$;
  \item every vertex $w\in C_j\setminus(\{v_j,v_{j-1}\})$, $j=2,\ldots,k-1$, should belong to the subinterval $[x_{j-1},\beta_{i-k+j-1})\sqcup\Delta_{i-k+j}\sqcup[\alpha_{i-k+j},x_{j})$;
  \item every vertex $w\in C_k\setminus\{v_{k-1}\}$ should belong to the subinterval $[x_{k-1},\beta_{i-1})\sqcup\Delta_{i}\sqcup\delta_i$.
\end{itemize}
The weights are independent, so denoting $y_j=x_j-\alpha_{i-k+j}$, we obtain the following estimate for the conditional probability:
\begin{align*}
&\left(\frac{1-p}{r}+y_1\right)^{n-1}\times\left(\frac{1-p}{r}+y_2+\frac{p}{r-1}-y_1\right)^{n-2}\times\ldots\times\\
&\times\left(\frac{1-p}{r}+y_{k-1}+\frac{p}{r-1}-y_{k-2}\right)^{n-2}\times\left(\frac{1-p}{r}+\frac{p}{r-1}-y_{k-1}+\frac{p}{r-1}\right)^{n-1}\leq\\
&\leq\text{|take out the factor $r^{-k(n-2)-2}$ and use the estimate $(1+y)^s\leq\exp\{ys\}$}|\leq\\
&\leq r^{-k(n-2)-2}e^{(n-1)\left(-p+ry_1\right)}e^{(n-2)\sum_{j=2}^{k-1}\left(-p+ry_j-ry_{j-1}+pr/(r-1)\right)}
e^{(n-1)\left(-p+pr/(r-1)-ry_{k-1}+\frac{p}{r-1}\right)}\leq\\
&+(n-2)\left(-p+ry_{k-2}-ry_{k-3}+pr/(r-1)\right)+(n-1)\left(-p+2pr/(r-1)-ry_{k-1}\right)\}=\\
&=r^{-k(n-2)-2}e^{(n-2)\left(-pk+\frac{kpr}{r-1}\right)}\cdot e^{\frac{2pr}{r-1}-2p}\cdot e^{ry_1-ry_{k-1}}\leq\\
&\leq |\text{since }0\leq y_1\leq p/(r-1)|\leq\\
&\leq r^{-k(n-2)-2}e^{(n-2)\left(-pk+\frac{kpr}{r-1}\right)}\cdot e^{\frac{3p}{r-1}}=\\
&= r^{-k(n-2)-2}e^{\frac {npk}{r-1}-\frac{2pk}{r-1}}\cdot e^{\frac{3p}{r-1}}\leq r^{-k(n-2)-2}e^{\frac {npk}{r-1}}\cdot e^{\frac{p}{r-1}}.
\end{align*}
To obtain the final estimate, we have to integrate over $(y_1,\ldots\,y_{k-1})\in [0,p/(r-1)]^{k-1}$ and substitute $p$ from \eqref{choice_p}. Thus, the probability under the consideration does not exceed
\begin{align*}
\left(\frac{p}{r-1}\right)^{k-1}& r^{-k(n-2)-2}e^{\frac {npk}{r-1}}\cdot e^{\frac{p}{r-1}}\leq\\
&\leq\left|\text{since $e^{p/(r-1)}\leq 2$}\right|\leq\\
&\leq r^{-k(n-1)-1}\left(\frac{\ln n}{n}\right)^{k-1}\cdot \left(\frac{n}{\ln n}\right)^{\frac kr}\cdot 2=\\
&= 2\left(\frac{\ln n}{n}\right)^{\frac{k(r-1)}{r}-1}r^{-(n-1)k-1}.
\end{align*}
So, the expected value of $X(i)$ can be estimated as follows:
$$
  {\sf E}X(i)\leq\sum_{k=1}^r 2{|E|\choose k}\cdot 2\left(\frac{\ln n}{n}\right)^{\frac{k(r-1)}{r}-1}\left(\frac 1r\right)^{(n-1)k+1}\leq\frac{4n}{r\ln n}\sum_{k=1}^\infty\frac{(0.01)^k}{k!}\leq \frac{0.04en}{r\ln n}.
	$$
The last statement of Lemma \ref{lemma4} obviously follows from Markov inequality:
$$
\Prob\left(\exists i:X(i)>\frac{n}{\ln n}\right)\leq r\frac{{\sf E}X(i)}{\frac{n}{\ln n}}\leq 0.04e.
$$
	
\end{proof1}

Let $K_1,\ldots,K_r$ denote the color classes of the initial coloring $C^0$. For every color class $K_{\alpha}$, let us define nonnegative integers $ex_{\alpha}$ and $sh_{\alpha}$, where $ex_{\alpha}$ stands for the excess value, i.e. a positive difference between number of vertices of the color ${\alpha}$ and $m/r$, and $sh_{\alpha}$ stands for the shortage value. Formally,
$$
  ex_{\alpha}=\begin{cases}|K_{\alpha}|-\frac mr,&\mbox{ if }|K_{\alpha}|-\frac mr>0;\\
  0,&\mbox{otherwise};
  \end{cases}
$$
$$
  sh_{\alpha}=\begin{cases}\frac mr-|K_{\alpha}|,&\mbox{ if }\frac mr-|K_{\alpha}|>0;\\
  0,&\mbox{otherwise}.
  \end{cases}
$$

\bigskip
The following lemma estimates the value of excesses and shortages.
\begin{lemma}\label{lemma5}
With probability at least $1/2-0.04e$ the following conditions hold simultaneously:
\begin{enumerate}
\item the first $r-1$ colors are in excess, i.e. the sizes if their color classes are at least $m/r$,
\item  for each color $i$  excess value $ex_i\leq\frac{p}{r(r-1)}m+\sqrt{\frac{13m\ln r}{r}}+\frac{n}{\ln n}$.

\end{enumerate}
\end{lemma}

\begin{proof1}
    Let $Z(i)$ denote the number of vertices, which belong to subinterval $\Delta_{i-1}\cup\delta_i$. Any such random variable has the binomial distribution $Bin(m,\frac{1-p}{r}+\frac{p}{r-1})$ with mean $m\left(\frac{1}{r}+\frac{p}{r(r-1)}\right)$. Then the number of vertices colored with color $i$ in $C^0$ is equal to
    $$
    Z(i)-X(i)+X(i-1).
    $$
	We use Chernoff's inequality which states that for any random binomial variable $\xi$ and any $z>0$, it holds that
\begin{align*}
	&\Prob(\xi < {\sf E}\xi - z)\leq\exp\left(-\frac{z^2}{2{\sf E}\xi}\right).\\
	&\Prob(\xi > {\sf E}\xi + z)\leq\exp\left(-\frac{z^2}{2({\sf E}\xi +z/3)}\right).
\end{align*}
Let us apply it for $\xi=Z(i)$ and $z =\sqrt{\frac{13m\ln r}{r}}$. Due to the initial restrictions on the parameter $r$ we have that $z<m/(2r)$. Hence,

\begin{align*}
&\Prob\left(\exists i:  Z(i) >m\left(\frac{1-p}{r}+\frac{p}{r-1}\right) + z\right)\leq r\exp\left(-\frac{z^2}{2(2\frac{m}{r}+z/3)}\right)<\\
&<r\exp\left(-\frac{z^2}{2\cdot\frac{13m}{6r}}\right)=r\exp\left(-3\ln r\right)=r^{-2}\leq 2^{-2}.
\end{align*}
Similarly,
$$
\Prob\left(\exists i: Z(i)<m\left(\frac{1-p}{r}+\frac{p}{r-1}\right)- z\right)\leq r^{-2}\leq 2^{-2}.
$$
Then, with probability at least $1-0.04e-2^{-1}=1/2-0.04e$ the following conditions hold simultaneously:
\begin{align*}
  &\frac{p}{r(r-1)}m+\sqrt{\frac{13m\ln r}{r}}+\frac{n}{\ln n}\geq Z(i)-X(i)+X(i-1)\\
  &Z(i)-X(i)+X(i-1)\geq m\left(\frac{1-p}{r}+\frac{p}{r-1}\right)-\sqrt{\frac{13m\ln r}{r}}-\frac{n}{\ln n}=\\
  &=\frac{m}{r}+\frac{mp}{r(r-1)}-\sqrt{\frac{13m\ln r}{r}}-\frac{n}{\ln n}.
\end{align*}
Hence, we have a desired upper bound on the excess value. Now notice that the first $(r-1)$ colors will be in excess if
$$
  \frac{mp}{r(r-1)}-\sqrt{\frac{13m\ln r}{r}}-\frac{n}{\ln n}=\frac{m\ln\left(\frac{n}{\ln n}\right)}{nr^2}-\sqrt{\frac{13m\ln r}{r}}-\frac{n}{\ln n}>0
$$
Since $m>\frac {n^2(r-1)}{2\ln n}$, the last inequality holds for $r<\sqrt[3]{\ln n}$ and all large enough $n$.
\end{proof1}

\section{Construction of an equitable coloring}

The final part of our work is devoted to restoring the color balance in the coloring $C^{0}$. According to Lemma \ref{lemma5} we know that with positive probability first $(r-1)$ color classes are in excess, i.e. their sizes are at least $m/r$. Therefore, we are going to recolor some vertices from the color classes $ K_{1}, ..., K_{r-1}$ with color $r$ keeping the lack of the monochromatic edges.

Now we formalize our idea: we want show that in every $\Delta_i$, $i=1,\ldots,r-1$, it is possible to choose a vertex subset $W_{i}$ of size $ex_{i}$ such that the recoloring of sets $W_{i}$ with color $r$ does not create monochromatic edges. If such choice of subsets is possible then we will be able to correct the color balance.

\subsection{Proof of the existence of sets $W_{i}$}

 Before we introduce the proof, let us simplify the next sum by one symbol and introduce a new parameter:
\begin{equation}\label{choice_q}
  q=\frac{mp}{r(r-1)}+2\sqrt{\frac{13m\ln r}{r}}+\frac{r+1}{r}\frac{n}{\ln n}
\end{equation}

 Now we are ready to establish the existence (with positive probability) of the required collection of sets. We consider a random subset $V_{i}$, $i=1\ldots,r-1$, constructed according to the binomial scheme: with probability
\begin{equation}\label{choice_tilde_p}
  \tilde{p}=q\times\frac{r}{(1-p)m}
\end{equation}
every vertex of $\Delta_i$ is included to $V_{i}$ independently of each other.

\bigskip
Which vertices of $V_i$ are not suitable to include in $W_{i}$? Clearly, we do not want to take the vertices whose recoloring with color $r$ can create monochromatic edges of color $r$. This situation happens when there is an edge $A$ and a subset of vertices $U(A)\subset A$ such that
\begin{itemize}
  \item all vertices in $A\setminus U(A)$ are colored with $r$ in the coloring $C^0$;
  \item all vertices in $U(A)$ belong to $\cup_{i=1}^{r-1}V_{i}$.
\end{itemize}
In this case  we do not want to take $U(A)$ completely into $\cup_{i=1}^{r-1}W_{i}$. Such an edge $A$ is called \emph{dangerous}. Let $Q$ denote the number of dangerous edges. So, to establish the existence of the required sets $W_1,\ldots,W_{r-1}$ it suffices to show that with positive probability the following conditions hold: for any $i=1,\ldots,r-1$,
$$
	|V_i|\geq Q+ex_{i}.
$$
We have already shown that with positive probability  $ex_{i}\leq\frac{p}{r(r-1)}m+\sqrt{\frac{13m\ln r}{r}}+\frac{n}{\ln n}$. So, it remains only to estimate $Q$ and the cardinalities $|V_i|$.

\bigskip
Let us start with $|V_i|$. It is clear from the construction that $|V_i|$ is a binomial random variable $Bin(m,q/m)$. So, by using Chernoff inequality, we get
\begin{align*}
   {\sf P}\left(|V_i|<q-\sqrt{\frac{13m\ln r}{r}}\right)&\le e^{-\frac{\frac{13m\ln r}{r}}{2q}}\le \\
   &\le|\mbox{since }q<\frac{2mp}{r(r-1)}|\le\\
   &\le e^{-(13)/4\cdot (r-1)(\ln r)p^{-1}}\le r^{-p^{-1}}<r^{-10}
\end{align*}
for all large enough $n$. Hence, with probability at least $1-r^{-9}\ge 1-2^{-9}$ every $|V_i|$ is at least
$$
  |V_i|\ge q-\sqrt{\frac{13m\ln r}{r}}=\frac{mp}{r(r-1)}+\sqrt{\frac{13m\ln r}{r}}+\frac{r+1}{r}\frac{n}{\ln n}.
$$

\subsection{Dangerous edges and $k$-complex ordered chains}

Suppose that $A$ is a dangerous edge. Let us denote $A'=A\setminus U(A)$. Then $A'$ can be considered as a monochromatic ``pseudo-edge'' of color $r$ during the evaluation of Algorithm 1. Hence, we can construct an ordered $k$-chain $H'=(C_1,...,C_k=A')$ for color $r$ with the ``pseudo-edge'' $A'$ as the last edge in the chain. Note that $A'$ should be contained in $\delta_{r-1}\cup \Delta_r$ ($A'$ is colored with $r$ in $C^0$) and, vice versa, none of the vertices of $U(A)$ can be contained in $\delta_{r-1}\cup \Delta_r$, so $A'=A\cap(\delta_{r-1}\cup \Delta_r)$. Now we can described a $k$-complex ordered chain $H''=(C_1,...,C_k=A)$ as follows:
\begin{itemize}
  \item $H'=(C_1,...,C_{k-1},A')$, where $A'=A\cap(\delta_{r-1}\cup\Delta_{r})$, is an ordered $k$-chain for color $r$ with pseudo-edge $A'$;
  \item all vertices in $U(A)=A\setminus A'$ belong to $\cup_{i=1}^{r-1}V_{i}$.
  \item every vertex $w\in A\cap C_j$ should belong to $V_{r-k+j}$,  $j\in\{1,...,k-2\}$;
  \item every vertex $w\in (A\setminus A')\cap C_{k-1}$ should belong to $V_{r-1}$.
\end{itemize}

\begin{lemma}\label{lemma6}
	For given edge $A$, the number of configurations that can form a $k$-complex ordered chain $k\geq 2$ in the hypergraph $H$ with $A$ as the last edge is at most $2\binom{|E|}{k-1}$.
\end{lemma}

\begin{proof1}
	Let us take an arbitrary unordered family of $k-1$ edges of $H$, say, $A_1,\ldots,A_{k-1}$. This can be done in at most $\binom{|E|}{k-1}$ ways. If $A_1,\ldots,A_{k-1}$ can form a chain then there are only two candidates for the the first edge in the chain (it should have only one common vertex with the union of others). After choosing the first edge the order in the chain is defined uniquely. \end{proof1}

\begin{lemma}\label{lemma7}
With probability at least $(1-0.02)$ the number of dangerous edges $Q$ does not exceed $\frac{1}{r}\left(\frac{n}{\ln n}\right)$
\end{lemma}

\begin{proof1}
Suppose $k\geq 2$. We want to estimate the probability that an ordered $k$-tuple of edges $H''=(C_1,\ldots,C_k=A)$ is a $k$-complex ordered chain. For any $j=1,\ldots,k-2$, we denote
$$
  v_j=C_j\cap C_{j+1} \mbox{ and } s_j=|A\cap C_j|.
$$
Then we choose a vertex $v_j\in C_{k-1}\cap A$ to be the unique common vertex of $A'$ and $C_{k-1}$ and denote $s_{k-1}=|A\cap C_{k-1}|-1$. Denote also
$$
  s=s_1+...+s_{k-1}.
$$

\bigskip
Recall that we use the notation $\delta_j=[\alpha_j,\beta_j)$. Thus, for given weights $\sigma(v_j)=x_j\in\delta_{r-k+j}$, $j=1,\ldots,k-1$, the event that $H''$ forms a $k$-complex ordered chain for color $r$, can be described as follows:
\begin{itemize}
  \item every vertex $w\in C_1\setminus\{v_1\}$ should belong to the subinterval $\Delta_{r-k+1}\sqcup[\alpha_{r-k+1},x_{1})$;
  \item every vertex $w\in C_j\setminus(\{v_j\cup v_{j-1}\})$, $j=2,\ldots,k-1$, should belong to the subinterval $[x_{j-1},\beta_{r-k+j-1})\sqcup\Delta_{r-k+j}\sqcup[\alpha_{r-k+j},x_{j})$;
  \item every vertex $w\in A\setminus(\{v_{k-1}\}\cup C_1\cup\ldots\cup C_{k-1})$ either belongs to the subinterval $[x_{k-1},\beta_{r-1})\sqcup\Delta_{r}$ or belongs to some $V_i$, $i=1,\ldots,r-1$;
  \item every vertex $w\in (A\cap C_j)$ should belong to $V_{r-k+j}$,  $j\in\{1,...,k-2\}$.
  \item every vertex $w\in (A\cap C_{k-1})\setminus \{v_{k-1}\}$ should belong to $V_{r-1}$.
\end{itemize}
Note that for every vertex $w$, ${\sf P}(w\in V_i)=(1-p)/r\cdot \widetilde{p}=q/m$ and this event is independent of the weights of other vertices.

\bigskip
As before, let $y_j=x_j-\alpha_{i-k+j}$. Now we are ready to estimate the conditional probability of event the $k$-complex chain occur given weights $x_j$ are fixed:

\begin{align*}
&\left(\frac{q}{m}\right)^s\left(\frac{1-p}{r}+y_1\right)^{n-1-s_1}\times\left(\frac{1-p}{r}+y_2+\frac{p}{r-1}-y_1\right)^{n-2-s_2}\times\ldots\times\\
&\times\left(\frac{1-p}{r}+y_{k-1}+\frac{p}{r-1}-y_{k-2}\right)^{n-2-s_{k-1}}
\times\left(\frac{1-p}{r}+\frac{p}{r-1}-y_{k-1}+\frac{(r-1)q}{m}\right)^{n-1-s}\leq
\end{align*}
(take out the factor $r$)
\begin{align*}
&\leq r^{-k(n-2)-2+s}\left(\frac{qr}{m}\right)^s
(1-p+ry_1)^{n-1-s_1}\left(1+\frac{p}{r-1}+ry_2-ry_1\right)^{n-2-s_2}\times\ldots\times\\
&\times\left(1+\frac{p}{r-1}+ry_{k-1}-ry_{k-2}\right)^{n-2-s_{k-1}}
\left(1+\frac{p}{r-1}-ry_{k-1}+\frac{(r-1)rq}{m}\right)^{n-1-s}\leq\\
&\leq \mbox{|since any expression in the brackets among the last three is at least $1-p$|}\leq\\
&\leq r^{-k(n-2)-2+s}\left(\frac{qr}{m}\right)^s(1-p+ry_1)^{n-1}\times\\
&\times\left(1+\frac{p}{r-1}+ry_2-ry_1\right)^{n-2}\times\ldots\times\left(1+\frac{p}{r-1}+ry_{k-1}-ry_{k-2}\right)^{n-2}\times\\
&\times\ldots\times\left(1+\frac{p}{r-1}-ry_{k-1}+\frac{r(r-1)q}{m}\right)^{n-1}\times\frac{1}{(1-p)^{s_1+...s_{k-1}+s}}\leq\\
&\leq r^{-k(n-2)-2+s}\left(\frac{qr}{m(1-p)^2}\right)^s\times\\
&\times\exp\left\{-p(n-1)+\frac{p(n-1)}{r-1}(k-1)+ry_1+\frac{r(r-1)(n-1)q}{m}\right\}.
\end{align*}
Let us analyze the obtained expression. Note that $q\le 2mp/(r-1)r$ and $(1-p)^2>1/2$, so \eqref{choice_p} implies that
$$
   \left(\frac{qr}{m(1-p)^2}\right)^s\leqslant \left(\frac{4p}{r-1}\right)^s\leqslant \left(\frac {4\ln n}{nr}\right)^s.
$$
Then since $y_1\leqslant \frac p{r-1}$
$$
  -p(n-1)+\frac{p(n-1)}{r-1}(k-1)+ry_1\leqslant 3p-pn+\frac{pn}{r-1}(k-1)=3p+\frac{k-r}{r-1}\ln\left(\frac{n}{\ln n}\right)^{\frac{r-1}{r}}.
$$
Finally, \eqref{choice_q} yields that
\begin{align*}
  &\frac{r(r-1)(n-1)q}{m}\leqslant \left(\frac{mp}{r(r-1)}+2\sqrt{\frac{13m\ln r}{r}}+\frac{r+1}{r}\frac{n}{\ln n}\right)\times\frac{r(r-1)(n-1)}{m}=\\
  &=p(n-1)+\frac{(r^2-1)n(n-1)}{m\ln n}+2\sqrt{\frac {13r\ln r}{m}}r(n-1)\le\\
  &\le|\mbox{using \eqref{vert_bound_new}}|\le p(n-1)+2(r+1)+2\sqrt{26\ln r\ln n}\cdot r\le \\
  &\le p(n-1)+\frac 1r\ln\frac{n}{\ln n}
\end{align*}
for $r<(\ln n)^{1/5}$ and all large enough $n$. Since $pn=\frac{r-1}{r}\ln\frac{n}{\ln n}$, we obtain the following estimate for the conditional probability:
\begin{align*}
  &r^{-k(n-2)-2+s} \left(\frac {4\ln n}{nr}\right)^se^{3p}\left(\frac{n}{\ln n}\right)^{\frac {k-r}{r}+1}=
  r^{-k(n-2)-2} \left(\frac {4\ln n}{n}\right)^se^{3p}\left(\frac{n}{\ln n}\right)^{\frac {k}{r}}.
\end{align*}

To obtain the final estimate, we have to integrate over the weights $y_1,\ldots\,y_{k-1}$ (factor $\left(p/(r-1)\right)^{k-1}$) and
sum up over all possible variants for the vertex $v_{k-1}$ (at most $s+1$ ways). Thus, we get the bound
\begin{align*}
  &r^{-k(n-2)-2} \left(\frac {4\ln n}{n}\right)^se^{3p}\left(\frac{n}{\ln n}\right)^{\frac {k}{r}}\left(\frac{p}{r-1}\right)^{k-1}(s+1)\le\\
  &\le(s+1)\left(\frac {4\ln n}{n}\right)^s r^{-k(n-1)-1}e^{3p}\left(\frac{n}{\ln n}\right)^{\frac {k}{r}-k+1}.
\end{align*}
Clearly, the bound is maximized for $s=1$.

If $k=1$ then the expected value of the number of dangerous edges can be estimated without complex ordered chains.
$$
\left(\frac{1-p}{r}+\frac{(r-1)q}{m}\right)^{n}\leq r^{-n}e^{-pn+\frac{(r-1)qrn}{m}}\le r^{-n}e^{\frac 1r\ln\frac{n}{\ln n}}.
$$
So, we are ready to estimate the expected value of the number of dangerous edges. Lemma \ref{lemma6} and the initial condition on the number edges imply that
\begin{align*}
  {\sf E}Q\le &|E| r^{-n}e^{\frac 1r\ln\frac{n}{\ln n}}+|E|\sum_{k=2}^r2{|E|\choose k-1}r^{-k(n-1)-1}e^{3p}\left(\frac{n}{\ln n}\right)^{\frac {k}{r}-k+1}\le\\
  &\le 0.01\cdot r^{-1}\frac {n}{\ln n}+2e^{3p}r^{-1}\frac{n}{\ln n}\sum_{k=2}^r\frac{(0.01)^k}{(k-1)!}\le\\
  &\le \frac {n}{r\ln n}\left(0.01+0.02e^{3p}(e^{0.01}-1)\right)\le 0.02\cdot\frac {n}{r\ln n}\\
\end{align*}
By Markov inequality we can conclude that the number of dangerous edges $Q$ does not exceed $\frac{1}{r}\left(\frac{n}{\ln n}\right)$ with probability at least $1-0.02$.
\end{proof1}

\subsection{Completion of the proof of Theorem \ref{thm:main}}

Let us sum up.

\bigskip
1)We have shown that the probability that there are monochromatic edges in the coloring $C^0$, does not exceed $0.04e$ (Lemma \ref{lemma_prob3}).

\bigskip
2) With probability at least $1/2-0.04e$ the first $r-1$ colors are in excess and for each color $i$, the excess value $ex_i\leq\frac{p}{r(r-1)}m+\sqrt{\frac{13m\ln r}{r}}+\frac{n}{\ln n}$. (Lemma \ref{lemma5}).

\bigskip
3) The probability that the number of dangerous edges is at least $\frac{1}{r}\left(\frac{n}{\ln n}\right)$ does not exceed $0.02$ (Lemma \ref{lemma7}).

4) With probability at least $1-2^{-9}$ every $V_i$ has cardinality at least $\frac{mp}{r(r-1)}+\sqrt{\frac{13m\ln r}{r}}+\frac{r+1}{r}\left(\frac{n}{\ln n}\right)$.

\bigskip
So, with probability at least $1/2-(0.04e+0.04e+0.02+2^{-9})>0$ there are no monochromatic edges and for every color $i$, the following relation holds:
$$
ex_{i}+Q\leq |V_{i}|.
$$
Hence, we can choose the required sets $W_{i}$ of size $ex_i$ and safely recolor them with color $r$ to obtain an equitable coloring of $H$. Theorem \ref{thm:main} is proved.

\section{Acknowledgements}

The work is supported by the grant of the President of Russian Federation no. MD-757.2019.1.

\renewcommand{\refname}{References}

\end{document}